\documentclass[11pt]{article}
\usepackage{amssymb,amsmath,amsthm,amsfonts}

\newcommand{\apos}{\Phi^+(A_{n-1})}
\newcommand{\bpos}{\Phi^+(B_n)}
\newcommand{\dpos}{\Phi^+(D_n)}
\newcommand{\bpone}{\Phi^+_1(B_n)}
\newcommand{\dpone}{\Phi^+_1(D_n)} 
\newcommand{\bptwo}{\Phi^+_2(B_n)}
\newcommand{\dptwo}{\Phi^+_2(D_n)}  
\newcommand{\dee}{\mathcal{D}}                                      

\newcommand{\gee}{\mathfrak{g}}

\newcommand{\hee}{\mathfrak{h}}

\newcommand{\C}{\mathbb{C}}

\newcommand{\Cash}{\mathbb{C}^{\#}}

\newcommand{\so}{\mathfrak{so}}

\newcommand{\en}{\mathfrak{n}}

\newtheorem{thm}{Theorem}[section]
\newtheorem{pro}[thm]{Proposition}
\newtheorem{lem}[thm]{Lemma}
\newtheorem{cor}[thm]{Corollary}

\newtheorem{defn}[thm]{Definition}

\begin{document}
 \title{Coadjoint orbits for $A_{n-1}^+, B_{n}^{+}$, and $D_{n}^{+}$}   
 \author{Shantala Mukherjee
 \thanks{This work is part of the author's doctoral dissertation, written
at the University of Wisconsin-Madison under the supervision of Prof. 
Georgia Benkart, and financially
supported in part by NSF grant \#{}DMS-0245082}\\
 Dept. of Mathematics \\
 DePaul University \\
 Chicago, IL 60614 \\
smukher2@condor.depaul.edu}

 \date{}
 \maketitle
 
 \begin{abstract}
 A complete description of the coadjoint orbits for $A_{n-1}^{+}$, the nilpotent Lie
 algebra of $n \times n$ strictly upper triangular matrices, has not yet been obtained, 
 though there has been steady progress on it ever since the orbit method was devised.
 We apply methods developed by Andr\'{e} to find defining equations for the 
 elementary coadjoint orbits for the maximal nilpotent Lie subalgebras of the orthogonal
 Lie algebras, and we also determine all the possible dimensions of coadjoint orbits in
 the case of $A_{n-1}^{+}$.

 \bigskip
 \noindent \textbf{MSC2000:} 17B30,17B35 
 \end{abstract}
 
 \section{Introduction} 
  The orbit method was created by Kirillov in an attempt to describe the unitary dual $\hat{N}_n$ for the nilpotent Lie group $N_n$ of 
$n \times n$ upper triangular matrices with $1$'s on the diagonal (the unitriangular group).
It turned out that the orbit method had much wider applications. In Kirillov's words:\\
`\dots all main questions of representation theory of Lie groups: construction of irreducible
representations, restriction-induction functors, generalized and infinitesimal characters, Plancherel
measure, etc., admit a transparent description in terms of coadjoint orbits' ({\cite{Kirillov:orb03}}).

The Lie algebra of $N_n$ is $\en_n$, which
consists of all $n \times n$ strictly upper triangular matrices. The group $N_n$ acts on the dual space $\en_n^*$ by the coadjoint
action, which will be explained later. A complete description of the set of coadjoint orbits $\en_n^*/N_n \simeq \hat{N}_n$
 in general is still not available, though progress had been made in the case of
the unitriangular group over a finite field. Andr\'{e} in {\cite{Andre:ja95}} defined ``basic sums" of 
elementary coadjoint orbits and their defining equations for the unitriangular group over an arbitrary
field and showed that the dual space $\en_n^*$
is a disjoint union of these basic sums of orbits. Similar results have been obtained by N.Yan, in his work on double orbits and cluster modules ({\cite{Yan:thesis01}}) of the unitriangular group over a finite field. 
Using Andr\'{e}'s results, we determine all possible dimensions of the coadjoint orbits of the
unitriangular group over $\C$. Elementary
coadjoint orbits (defined later) are in some sense the ``smallest" coadjoint orbits. Adapting the
methods of Andr\'{e}, we derive the defining equations for the elementary coadjoint orbits in the case
that the Lie algebra is a maximal nilpotent subalgebra of an orthogonal Lie algebra
over the field $\C$ of complex numbers. The case of the symplectic Lie algebra has so
far failed to yield a consistent pattern, therefore it is not discussed here.

\section{The Lie algebras $A_{n-1}^{+},\ B_{n}^{+},\ D_{n}^{+}$}
                                                                        

Let $L =
\mathfrak{sl}_{n}(\mathbb{C})$, the Lie algebra of $n \times n$
complex matrices of trace $0$. Let $E_{ij}$ denote the standard matrix 
unit with $1$ in the $(i,j)$ position and $0$ elsewhere. Thus $L$
has a basis 
\[ \{E_{ij} \mid 1\leq i\neq j\leq n\} \cup \{E_{ii} - E_{i+1,i+1}
\mid i = 1,\dots,n-1\}.\]

 Relative to the Cartan subalgebra
$\mathfrak{h}$ spanned by the diagonal matrices $E_{ii}-E_{i+1,i+1},\
i=1,\dots,n-1$; $L$ decomposes into root spaces (common eigenspaces).
Thus
\[ L = \hee \oplus \bigoplus_{1\leq i\neq j\leq
n}L_{\epsilon_i-\epsilon_j},\]
where $\epsilon_i \colon \hee \rightarrow \C$ denotes the projection 
onto the $(i,i)$ entry, and 
\[ L_{\epsilon_i-\epsilon_j} = \{x \in L \mid
[h,x]=(\epsilon_i-\epsilon_j)(h)x \ \forall\ h \in \hee\} = \C E_{ij}. \]
The roots $\epsilon_i-\epsilon_j$ are linear combinations of  
the simple roots 
$\epsilon_1 - \epsilon_2, \epsilon_2 - \epsilon_3,\dots,\epsilon_{n-1} - 
\epsilon_n$, and the coefficients are either all nonpositive or all
nonnegative integers.
The positive roots are given by 
\[ \apos = \{\epsilon_i - \epsilon_j \mid 1 \leq i < j \leq n \}. \]
Thus, there are $\frac{1}{2}(n-1)n$ positive roots. 

The sum $\bigoplus_{\alpha \in \apos} L_{\alpha}$ of the root spaces
corresponding to the positive roots is the 
nilpotent Lie algebra $\en_n$ of strictly upper triangular matrices. Here we denote
this Lie algebra by $A_{n-1}^+$. It has a basis of root vectors $\{e_{\alpha} \mid \alpha \in \apos\}$
where $e_{\alpha} = E_{ij}$ for $\alpha = \epsilon_i-\epsilon_j,\ 1\leq i<j\leq n$.

\bigskip

Now we establish our conventions for  the root systems $B_n$ and $D_n$. The following definitions can be found
in {\cite[Sec.~18.1]{Ful:reptheory91}}. Let $J_m$ be the $m \times m$ matrix with $1$'s along the antidiagonal
and $0$'s elsewhere.  The orthogonal Lie algebra $\mathfrak{so}_m(\C)$ is the Lie algebra of
$m \times m$ matrices $X$ satisfying the relation $X^tJ_m + J_mX = 0$. Thus the matrices
in $\mathfrak{so}_m(\C)$ are antisymmetric about the antidiagonal.

Relative to the Cartan subalgebra
$\mathfrak{h}$ spanned by the diagonal matrices $E_{ii} - 
E_{2n+2-i,2n+2-i}$, 
the odd orthogonal Lie algebra
$L=\mathfrak{so}_{2n + 1}(\mathbb{C})$ decomposes into root spaces
\[ L = \hee \oplus \bigoplus_{1\leq i\neq j\leq n}L_{\epsilon_i \pm \epsilon_j}\oplus \bigoplus_{1\leq i\leq n}L_{\epsilon_i},\] 
where $\epsilon_i \colon \hee \rightarrow \C$ denotes the projection onto the
$(i,i)$ entry. The roots are linear combinations of the 
simple roots 
$\epsilon_1 - \epsilon_2, \epsilon_2 - \epsilon_3,\dots,\epsilon_{n-1} - 
\epsilon_n, \epsilon_n$, with coefficients that are either all nonnegative or all nonpositive integers.
The positive roots are given by
\[ \bpos = \{\epsilon_i \pm \epsilon_j 
\mid 1 \leq i < j \leq n \} \cup \{\epsilon_i \mid 1 \leq i \leq n \}  \]
\noindent and $\vert \bpos \vert = n^2$.

We partition the set of positive roots into two subsets
\begin{defn}
$\bpone  = \{\epsilon_i - \epsilon_j \mid 1\leq i<j\leq n\} \cup \{\epsilon_i \mid i = 1,\dots,n\}$
\end{defn}
\begin{defn}
$\bptwo  = \{\epsilon_i + \epsilon_j \mid 1 \leq i < j \leq n \}$
\end{defn}

The sum of spaces $\bigoplus_{\alpha \in \bpos}L_{\alpha}$ is a 
finite-dimensional nilpotent Lie algebra consisting of all strictly upper
triangular matrices of size $2n + 1$ which are anti-symmetric about the
antidiagonal (and have zeroes on the antidiagonal). We denote this Lie algebra by $B_n^+$.
It has a basis of root vectors $\{e_{\alpha} \mid \alpha \in \bpos \}$ where 
\begin{equation}\label{E:Brootvectors} 
e_{\alpha} =
\begin{cases}
E_{ij}-E_{2n+2-j,2n+2-i},   &\text{if $\alpha = \epsilon_i - \epsilon_j,\ \ 1\leq i<j\leq n$;} \\
E_{i,n+1}-E_{n+1,2n+2-i},  &\text{if $\alpha = \epsilon_i,\ \ 1\leq i \leq n$;} \\
E_{i,2n+2-j}-E_{j,2n+2-i},  &\text{if $\alpha = \epsilon_i + \epsilon_j,\ \ 1\leq i< j\leq n$.}
\end{cases} \end{equation}

\bigskip

Relative to the Cartan subalgebra 
$\mathfrak{h}$ spanned by the matrices $H_i = E_{ii} - E_{2n+1-i,2n+1-i}$, the even orthogonal
Lie algebra $L=\so_{2n}(\C)$ decomposes into root spaces
\[ L = \hee \oplus \bigoplus_{1\leq i\neq j\leq n}L_{\epsilon_i\pm\epsilon_j} \]
where  $\epsilon_i \colon \hee \rightarrow \C$ denotes the projection onto the
$(i,i)$ entry. The roots are linear combinations of the 
simple roots 
$\epsilon_1 - \epsilon_2, \epsilon_2 - \epsilon_3,\dots,\epsilon_{n-1} - 
\epsilon_n, \epsilon_{n-1}+\epsilon_n$, and the coefficients are either all nonnegative or all nonpositive integers.
The positive roots are given by
\[ \dpos = \{\epsilon_i \pm \epsilon_j 
\mid 1 \leq i < j \leq n \}  \]
\noindent and $\vert \dpos \vert = n^2-n$.

The elements of $\dpos$ can be partitioned into two subsets:
\begin{defn}
$\dpone  = \{ \epsilon_i - \epsilon_j \mid 1\leq i<j\leq n \},$
\end{defn}
\begin{defn}
$\dptwo = \{ \epsilon_i + \epsilon_j \mid 1\leq i<j\leq n \}$
\end{defn}

The sum of root spaces $\bigoplus_{\alpha \in \dpos}L_{\alpha}$ is a      
finite-dimensional nilpotent Lie algebra consisting of all strictly upper
triangular matrices of size $2n$ which are anti-symmetric about the
antidiagonal (and have zeroes on the antidiagonal). We denote this Lie algebra by $D_n^+$.
It has a basis of root vectors $\{e_{\alpha} \mid \alpha \in \dpos \}$ where
\begin{equation}\label{E:Drootvectors} 
e_{\alpha} = 
\begin{cases}
E_{ij} - E_{2n+1-j,2n+1-i},  &\text{if $\alpha = \epsilon_i - \epsilon_j,\ 1\leq i < j\leq n$;} \\
E_{i,2n+1-j} - E_{j,2n+1-i},  &\text{if $\alpha = \epsilon_i + \epsilon_j,\ 1\leq i < j \leq n$.}
\end{cases} \end{equation}

\section{Singular and Regular Roots}

For a positive root $\alpha$ in any one of the sets $\apos, \bpos,$ or $\dpos$,
we will define two sets: $S(\alpha)$, the set of
$\alpha$-\emph{singular roots}; and $R(\alpha)$, the set of
$\alpha$-\emph{regular roots}. The set $S(\alpha)$ is the union of all pairs of positive roots which
sum up to $\alpha$.

For a positive root of the form $\epsilon_i - \epsilon_j,\ 1\leq i< j\leq n$, in the sets $\apos$, $\bpone$ , and $\dpone$, we have
\[ S(\epsilon_i - \epsilon_j) = 
\begin{cases}
\bigcup_{k=i+1}^{j-1}\{\epsilon_i - \epsilon_k, \epsilon_k - \epsilon_j\},  &\text{if $ j-i>1$} \\
\emptyset,  &\text{otherwise}
\end{cases}
 \]
We see that $\vert S(\epsilon_i - \epsilon_j) \vert = 2(j - i-1)$.  

For a positive root of the form $\epsilon_i,\ 1\leq i\leq n$, in
$\bpone$ we have
\[ S(\epsilon_i) = 
\begin{cases}
\bigcup_{k=i+1}^{n}\{\epsilon_i - \epsilon_k,\epsilon_k\},  &\text{if
$1\leq i\leq n-1$} \\
\emptyset,  &\text{if $i = n.$}
\end{cases} \]
Here $\vert S(\epsilon_i) \vert = 2(n-i)$.
\bigskip

For the roots $\epsilon_i + \epsilon_j,\ 1\leq i<j\leq n$, in $\bptwo$ 
we have:
\begin{equation}
\begin{split}
S(\epsilon_i + \epsilon_j)  = & 
\bigcup_{k=i+1}^{j-1}\{\epsilon_i - \epsilon_k, \epsilon_k + \epsilon_j\} 
\cup \bigcup_{k=j+1}^{n}\{\epsilon_i - \epsilon_k, \epsilon_j + \epsilon_k\} \\ 
& \cup \{\epsilon_i, \epsilon_j\}
\cup \bigcup_{k=j+1}^{n}\{\epsilon_i + \epsilon_k, \epsilon_j - \epsilon_k\} 
\end{split}
\end{equation}
Notice that $\vert S(\epsilon_i + \epsilon_j) \vert =  2(2n-(i+j))$.


For the positive roots $\epsilon_i + \epsilon_j,\ 1\leq i<j \leq n$, in 
$\dptwo$ we have:
\begin{equation}
\begin{split}
S(\epsilon_i + \epsilon_j)  = & 
\bigcup_{k=i+1}^{j-1}\{\epsilon_i - \epsilon_k, \epsilon_k + \epsilon_j\} 
\cup \bigcup_{k=j+1}^{n}\{\epsilon_i - \epsilon_k, \epsilon_j + \epsilon_k\} \\ 
&\cup \bigcup_{k=j+1}^{n}\{\epsilon_i + \epsilon_k, \epsilon_j - \epsilon_k\} 
\end{split}
\end{equation}
In this case, $\vert S(\epsilon_i+\epsilon_j)\vert = 2(2n-i-j-1)$.

For $\alpha \in \apos, \bpos,$ or $\dpos$, we define $R(\alpha)$ to be the complement
of the set $S(\alpha)$ in the respective set of positive roots. Clearly, $\alpha \in R(\alpha)$.

\section{Elementary Coadjoint Orbits}\label{S:elemcoadorb}
The definitions and results in this subsection are taken from {\cite[Sec.~1]{Andre:ja95}}.
Let $\Phi^+$ denote one of the sets of positive roots $\apos, \bpos$, or
$\dpos$. Let $\gee$ denote the corresponding nilpotent Lie algebra
$A_{n-1}^+, B_n^+$, or $D_n^+$.
The  group $G = \exp(\gee)$
acts on the dual space $\gee^*$ by the \emph{coadjoint action}
\[ (g.f)(x) = f(g^{-1}xg),\ \ \forall g \in G,\ f \in \gee^*,\ x \in \gee.\]
Then, by {\cite[Thm. 6.2.4]{Dixmier:ea96}}, there is a one-one
correspondence between the $G$-orbits in $\gee^*$ and the primitive
ideals of $U(\gee)$, which are the annihilators of the simple
$U(\gee)$-modules, constructed as follows: For any $f$ in the $G$-orbit, a simple $\gee$-module is
obtained by inducing a one-dimensional module of a maximal subalgebra of $\gee$ that is subordinate
to $f$, up to a $\gee$-module ({\cite[Thm.~6.1.1]{Dixmier:ea96}}). The annihilator in $U(\gee)$ of
this module is the primitive ideal $I(f)$ that corresponds to the
$G$-orbit $\Omega_f$ of $f$.
By {\cite[Thm. 4.7.9(iii)]{Dixmier:ea96}}, factoring $U(\gee)$ by a primitive
ideal gives a Weyl
algebra $\mathcal{A}_m$, which is the non-commutative algebra of
algebraic differential operators on 
a polynomial ring $\C[x_1,\dots,x_m]$.
\begin{defn}
$\gee^f = \{x \in \gee \mid f([x,y])=0\ \forall y \in \gee\}$
is the radical of the form $f$.
\end{defn}
Then $\dim \gee / \gee^f = \dim \Omega_f = 2m$, where $U(\gee)/I(f) = \mathcal{A}_m$. 

For any $\alpha \in \Phi^+$, let $e_{\alpha}^*$ denote the element of the
dual vector space $\gee^*$ defined as follows: for any $\beta \in \Phi^+$,

\[
  e_{\alpha}^*(e_{\beta}) = 
\begin{cases}
   1 ,      &\text{if $\alpha = \beta$;} \\
   0 ,      &\text{otherwise.}
\end{cases}
\]

\noindent Then $\{e_{\alpha}^* \mid \alpha \in \Phi^+ \}$ is a basis of
$\gee^*$.
Let $c \in \C$ be non-zero. Then, under
the coadjoint action of
the group $G$ ($=\exp(\gee)$) on $\gee^*$, the
coadjoint orbit $O_{\alpha}(c)$ that contains the element
$ce_{\alpha}^*$ 
is called
the $\alpha$-th \emph{elementary orbit associated with $c$}.
Note that if $f = ce_{\alpha}^*$, then $\gee^{f}=\{x \in \gee \mid
f([x,y])=0\ \forall \ y \in \gee\}$ is spanned by the 
$e_{\beta},\ \beta \in R(\alpha)$,
so 
\begin{equation}
\dim(O_{\alpha}(c)) = \dim(\gee/\gee^f) = \dim \gee - \dim \gee^f=\vert \Phi^+ \vert - \vert
R(\alpha) \vert =\vert S(\alpha) \vert. 
\end{equation}

Let $t$ be an arbitrary scalar in $\C$ and let $\beta \in \Phi^+$ . 
Then, the matrix $\exp(te_{\beta}) \in G$.
So, by the definition of the coadjoint representation, we have for
any $\gamma \in \Phi^+$:
\begin{equation}
\begin{split} 
\exp(te_{\beta}).e_{\alpha}^*(e_{\gamma}) & = 
e_{\alpha}^*(\exp(-te_{\beta})e_{\gamma}\exp(te_{\beta})) \\ 
& = e_{\alpha}^*(\exp(ad(-te_{\beta})(e_{\gamma})) \\  
& = e_{\alpha}^*(e_{\gamma} - t[e_{\beta},e_{\gamma}] + \frac{1}{2}t^2
[e_{\beta},[e_{\beta},e_{\gamma}]] + \dots ) .
\end{split}
\end{equation}

It is clear that for any simple root $\alpha_i \in \Phi^+$ and any 
$c \in \C$, the coadjoint orbit containing $ce_{\alpha_i}^*$ is
equal to $\{ce_{\alpha_i}^*\}$, because 
$ce_{\alpha_i}^*([\gee,\gee]) = 0.$

By Prop.~8.2 in {\cite{Hum:linalggps75}}, any coadjoint orbit is an irreducible variety
in $\gee^*$, so in particular, the elementary coadjoint orbit $O_{\alpha}(c)$ is an irreducible variety of dimension $\vert S(\alpha) \vert$.

\newpage

Andr\'{e} in {\cite[Lem.~2]{Andre:basic95}} describes the elementary orbit 
$O_{\alpha}(c)$, for any
$\alpha \in \apos$ and any non-zero scalar $c$.

If $g \in G$, then $g.(ce_{\alpha}^*) = c(g.e_{\alpha}^*)$, so
\[ f \in O_{\alpha}(c)\ \ \text{if and only if}\ \ \frac1c f \in 
O_{\alpha}(1). \]
Thus, it is enough to determine the defining equations for 
$O_{\alpha}(1)$.
Adapting Andr\'{e}'s proof, we obtain the  
defining equations for
elementary orbits $O_{\alpha}(1)$, where $\alpha \in \apos, \bpone$, or $\dpone$:

\begin{thm}
\begin{description}
\item[(a)]
Let $\alpha = \epsilon_i - \epsilon_j \in \apos, \bpone$, or $\dpone$, where $1\leq i<j\leq n$. Let $\gee$ denote
the corresponding nilpotent Lie algebra $A_{n-1}^+, B_n^+$, or $D_n^+$. Then $O_{\alpha}(1)$ consists of all elements 
$f \in \mathfrak{g}^*$ which satisfy the equations 
\begin{equation}\label{E:orbB}
 f(e_{\beta}) = 
   \begin{cases} 
    1,    &\text{if $\beta = \alpha$;}\\
    f(e_{\epsilon_i - \epsilon_s})f(e_{\epsilon_r - \epsilon_j}),
    &\text{if $\beta = \epsilon_r - \epsilon_s,\ \ i<r<s<j$;}\\
    0   &\text{otherwise,}
   \end{cases} 
\end{equation}
for $\beta \in R(\alpha)$, and $f$ takes arbitrary values on $e_{\beta}$ for $\beta \in S(\alpha)$. 
\item[(b)]
Let $\gee$ denote the nilpotent Lie algebra $B_n^+$. For $\epsilon_i \in \bpone,\ 1\leq i\leq n$,
the elementary orbit $O_{\epsilon_i}(1)$ consists of all elements $f \in \gee^*$ that satisfy the
equations
\begin{equation}\label{E:orbB1}
f(e_{\beta}) = 
\begin{cases}
1, &\text{if $\beta = \epsilon_i$}; \\
f(e_{\epsilon_i - \epsilon_s})f(e_{\epsilon_r}), &\text{if $\beta = \epsilon_r - \epsilon_s,\ i< r<s\leq n$};\\
0 &\text{otherwise},
\end{cases}
\end{equation}
for $\beta \in R(\epsilon_i)$, and $f$ takes arbitrary values on $e_{\beta}$ for $\beta \in S(\epsilon_i)$.
\end{description}
\end{thm} 

\bigskip

\bigskip

\noindent \textbf{Proof.} 
\textbf{(a)}
Let $\mathcal{V}$ be the variety in $\mathfrak{g}^*$ consisting of all $f \in \gee^*$ that satisfy
the equations~(\ref{E:orbB}) and let $f \in \mathcal{V}$. Then 
\[ f = \left (\prod_{k=i+1}^{j-1}\exp(f(e_{\epsilon_i - \epsilon_k})e_{\epsilon_k - \epsilon_j})\prod_{k= 
i+1}^{j-1}\exp(-f(e_{\epsilon_k - \epsilon_j})e_{\epsilon_i - \epsilon_k})\right ).e_{\alpha}^* \in 
O_{\alpha}(1) ,\]
so $\mathcal{V} \subseteq O_{\alpha}(1)$.
To show that  equality holds, let $T:\mathcal{V} \rightarrow \mathbb{C}^{2(j-i-1)}$
be the map defined by applying $f$ to pairs in $S(\alpha)$ as follows:
\[ T(f) = 
\left(f(e_{\epsilon_i - \epsilon_{i+1}}),f(e_{\epsilon_{i+1} -
    \epsilon_j}),\dots,
f(e_{\epsilon_i - \epsilon_{j-1}}),f(e_{\epsilon_{j-1}-\epsilon_j})\right )\]
for all $f \in \mathcal{V}$. 

\noindent (For example, if $\gee=A_3^+$ and $\alpha =
    \epsilon_1-\epsilon_4$, then 
\[T(f) = \left(f(e_{\epsilon_1-\epsilon_2}),f(e_{\epsilon_2-\epsilon_4}),f(e_{\epsilon_1-\epsilon_3}),
f(e_{\epsilon_3-\epsilon_4})\right).)\]
Then $T$ is an isomorphism of algebraic varieties, and since $\C ^{2(j-i-1)}$ is an irreducible variety, it follows that $\mathcal{V}$ is 
irreducible and $\dim\mathcal{V} = 2(j-i-1)$.

The coadjoint orbit $O_{\alpha}(1)$ is also an irreducible algebraic variety of 
dimension $2(j-i-1)$. We have two irreducible varieties $\mathcal{V}$ and $O_{\alpha}(1)$ of the same dimension and $\mathcal{V} \subseteq O_{\alpha}(1)$, so it follows that $\mathcal{V} = O_{\alpha}(1)$.

\textbf{(b)} Let $\mathcal{V}$ be the variety in $\mathfrak{g}^*$ consisting of all $f \in \gee^*$ that satisfy
the equations~(\ref{E:orbB1}). Let $f \in \mathcal{V}$. Then 
\[ f = \left (\prod_{k=i+1}^{n}\exp(f(e_{\epsilon_i - \epsilon_k})e_{\epsilon_k})\prod_{k= 
i+1}^{n}\exp(-f(e_{\epsilon_k})e_{\epsilon_i - \epsilon_k})\right ).e_{\alpha}^* \in 
O_{\alpha}(1) ,\]
so $\mathcal{V} \subseteq O_{\alpha}(1)$.
To show that  equality holds, let $T:\mathcal{V} \rightarrow \mathbf{C}^{2(n-i)}$
be the map defined by:
\[ T(f) = 
\left(f(e_{\epsilon_i - \epsilon_{i+1}}),f(e_{\epsilon_{i+1}}),
\dots,f(e_{\epsilon_i - \epsilon_n}),f(e_{\epsilon_n})\right )\]
for all $f \in \mathcal{V}$. 

\noindent (For e.g., if $\gee = B_3^+$ and $\alpha = \epsilon_1$, then
\[T(f) = \left
  (f(e_{\epsilon_1-\epsilon_2}),f(e_{\epsilon_2}),f(e_{\epsilon_1-\epsilon_3}),
f(e_{\epsilon_3})\right ).)\]
Then $T$ is an isomorphism of algebraic varieties, and because $\C ^{2(n-i)}$ is an irreducible variety, $\mathcal{V}$ is 
irreducible and $\dim\mathcal{V} = 2(n-i)$.

The coadjoint orbit $O_{\alpha}(1)$ is also an irreducible algebraic variety of 
dimension $2(n-i)$. We have two irreducible varieties $\mathcal{V}$ and $O_{\alpha}(1)$ of the same dimension and $\mathcal{V} \subseteq O_{\alpha}(1)$, so it follows that $\mathcal{V} = O_{\alpha}(1)$.
\qed

\bigskip

Next, we can describe the defining equations of the elementary orbit $O_{\alpha}(1)$ for the positive
roots
$\alpha = \epsilon_i + \epsilon_j,\ 1\leq i<j\leq n$ in $\bptwo$ or $\dptwo$.
\begin{thm}
\begin{description}
\item[(i)] Let $\gee = B_n^+$. For a positive root $ \epsilon_i + \epsilon_j,\ 1\leq i<j\leq n$ in $\bptwo$,
$f \in O_{\epsilon_i +\epsilon_j}(1)$ if and only if $f$ satisfies:
\begin{equation}\label{E:orbB2}
f(e_{\beta}) =
\begin{cases}
1\ \ \ \ \ \ \ \ \ \ \ \ \ \ \ \ \ \ \ \ \ \text{if $\beta = \epsilon_i + \epsilon_j$;} \\
f(e_{\epsilon_i-\epsilon_s})f(e_{\epsilon_r+\epsilon_j}) 
\ \ \text{if $\beta = \epsilon_r-\epsilon_s,\ i\leq r<s\leq j$;} \\
f(e_{\epsilon_r+\epsilon_j})\left (- \frac{1}{2}{f(e_{\epsilon_i})^2} + 
\sum_{k=j+1}^{n}(-1)^kf(e_{\epsilon_i-\epsilon_k})f(e_{\epsilon_i+\epsilon_k}) \right ) \\
\ \ \ \ \ \ \ \ \ \ \ \ \ \ \ \ \ \ \ \ \ \ \ \text{if $\beta = \epsilon_r - \epsilon_j,\ i\leq r<j$;} \\
f(e_{\epsilon_i \pm \epsilon_s})f(e_{\epsilon_r+\epsilon_j}) 
\ \ \text{if $\beta = \epsilon_r \pm \epsilon_s,\ i < r<j<s\leq n$;} \\
f(e_{\epsilon_j \pm \epsilon_s})f(e_{\epsilon_i+\epsilon_r}) - f(e_{\epsilon_i \pm \epsilon_s})f(e_{\epsilon_j + \epsilon_r}) \\
\ \ \ \ \ \ \ \ \ \ \ \ \ \ \ \ \ \ \ \ \ \ \ \text{if $\beta = \epsilon_r \pm \epsilon_s,\ j<r<s\leq n$;} \\
f(e_{\epsilon_i})f(e_{\epsilon_r + \epsilon_j}) 
\ \ \ \ \ \ \text{if $\beta = \epsilon_r,\ i < r < j$;} \\
f(e_{\epsilon_j})f(e_{\epsilon_i + \epsilon_r})-f(e_{\epsilon_i})f(e_{\epsilon_j+\epsilon_r})
\ \ \text{if $\beta = \epsilon_r,\ j<r\leq n$;} \\
0 
\ \ \ \ \ \ \ \ \ \ \ \ \ \ \ \ \ \ \ \ \ \ \text{otherwise,}
\end{cases}
\end{equation}
for all $\beta \in R(\epsilon_i + \epsilon_j)$, and $f$ takes arbitrary values on $e_{\beta}$ for all
$\beta \in S(\epsilon_i + \epsilon_j)$.

\item[(ii)]Let $\gee = D_n^+$. For a positive root $ \epsilon_i + \epsilon_j,\ 1\leq i<j\leq n$ in $\dptwo$,
$f \in O_{\epsilon_i +\epsilon_j}(1)$ if and only if $f$ satisfies:

\bigskip

\begin{equation}\label{E:orbD2}
f(e_{\beta}) =
\begin{cases}
1 \ \ \ \ \ \ \ \ \ \ \ \ \ \ \ \ \ \ \ \ \ \text{if $\beta = \epsilon_i + \epsilon_j$;} \\
f(e_{\epsilon_i-\epsilon_s})f(e_{\epsilon_r+\epsilon_j}) 
\ \ \text{if $\beta = \epsilon_r-\epsilon_s,\ i\leq r<s\leq j$;} \\
f(e_{\epsilon_r+\epsilon_j})\left (\sum_{k=j+1}^{n}(-1)^kf(e_{\epsilon_i-\epsilon_k})f(e_{\epsilon_i+\epsilon_k}) \right ) \\
\ \ \ \ \ \ \ \ \ \ \ \ \ \ \ \ \ \ \ \ \ \ \ \text{if $\beta = \epsilon_r - \epsilon_j,\ i\leq r<j$;} \\
f(e_{\epsilon_i \pm \epsilon_s})f(e_{\epsilon_r+\epsilon_j}) 
\ \ \text{if $\beta = \epsilon_r \pm \epsilon_s,\ i < r<j<s\leq n$;} \\
f(e_{\epsilon_j \pm \epsilon_s})f(e_{\epsilon_i+\epsilon_r}) - f(e_{\epsilon_i \pm \epsilon_s})f(e_{\epsilon_j + \epsilon_r}) \\
\ \ \ \ \ \ \ \ \ \ \ \ \ \ \ \ \ \ \ \ \ \ \ \text{if $\beta = \epsilon_r \pm \epsilon_s,\ j<r<s\leq n$;} \\
0 \ \ \ \ \ \ \ \ \ \ \ \ \ \ \ \ \ \ \ \ \ \ \text{otherwise,}
\end{cases}
\end{equation}
for all $\beta \in R(\epsilon_i + \epsilon_j)$, and $f$ takes arbitrary values on $e_{\beta}$ for all
$\beta \in S(\epsilon_i + \epsilon_j)$.
\end{description}
\end{thm}

\noindent \textbf{Proof.}
\textbf{(i)} Recall that
\begin{equation*}
\begin{split}
S(\epsilon_i + \epsilon_j)  = & 
\bigcup_{k=i+1}^{j-1}\{\epsilon_i - \epsilon_k, \epsilon_k + \epsilon_j\} 
\cup \bigcup_{k=j+1}^{n}\{\epsilon_i - \epsilon_k, \epsilon_j + \epsilon_k\} \\ 
& \cup \{\epsilon_i, \epsilon_j\}
\cup \bigcup_{k=j+1}^{n}\{\epsilon_i + \epsilon_k, \epsilon_j - \epsilon_k\} ,
\end{split}
\end{equation*}
and $\vert S(\epsilon_i + \epsilon_j) \vert =  2(2n-(i+j))=\dim O_{\epsilon_i+\epsilon_j}(1)$.
The set $S(\epsilon_i+\epsilon_j)$ can be written as the disjoint union of two subsets
$S_{(i)}$ and $S_{(j)}$ defined as follows:
\[ S_{(i)} = \bigcup_{k=i+1}^{j-1}\{\epsilon_i - \epsilon_k\} \cup \bigcup_{k=j+1}^n\{\epsilon_i-\epsilon_k\}\cup 
\{\epsilon_i\}\cup \bigcup_{k=j+1}^n\{\epsilon_i+\epsilon_k\}\] and
\[ S_{(j)} = \bigcup_{k=i+1}^{j-1}\{\epsilon_k+\epsilon_j\}\cup \bigcup_{k=j+1}^n\{\epsilon_j+\epsilon_k\}\cup 
\{\epsilon_j\}\cup \bigcup_{k=j+1}^n\{\epsilon_j-\epsilon_k\}.\]
For $\gamma \in S_{(i)}$, let $\gamma'$ be the unique element of $S_{(j)}$ such that
$\gamma + \gamma' = \epsilon_i + \epsilon_j$. Let $n(\gamma,\gamma')$ be the integer 
such that $[e_{\gamma},e_{\gamma'}] = n(\gamma,\gamma')e_{\epsilon_i+\epsilon_j}$
It follows from equations~(\ref{E:Brootvectors}) that $n(\gamma,\gamma') = \pm 1$.

Let $\mathcal{V}$ be the variety in $\gee^*$ consisting of all $f \in \gee^*$ 
satisfying equations~(\ref{E:orbB2}). Let $f \in \mathcal{V}$. Then
\[ f = \left ( \prod_{\gamma,\gamma' }\exp(n(\gamma,\gamma')f(e_{\gamma})e_{\gamma'}) 
\prod_{\gamma,\gamma' }\exp(-n(\gamma,\gamma')f(e_{\gamma'})e_{\gamma})\right ).e_{\epsilon_i+\epsilon_j}^* ,\]
where $\gamma \in S_{(i)}$ and $\gamma' \in S_{(j)}$.

Thus $\mathcal{V} \subseteq O_{\epsilon_i+\epsilon_j}(1)$. 
To show equality, let us define a map 
$T:\mathcal{V} \rightarrow \mathbb{C}^{2(2n-(i+j))}$ by 
\[ T(f) =
\left(f(e_{\gamma_1}),f(e_{\gamma_1'}),f(e_{\gamma_2}),f(e_{\gamma_2'}),
  \dots,
f(e_{\gamma_{2n-(i+j)}}),f(e_{\gamma_{2n-(i+j)}'})\right)\]
where $\gamma_k \in S_{(i)}$, for all $k$. 

\noindent (For example, if we have $\epsilon_1+\epsilon_3 \in
  \Phi_2^+(B_3)$, then 
\[T(f) = \left (f(e_{\epsilon_1-\epsilon_2}),f(e_{\epsilon_2+\epsilon_3}),f(e_{\epsilon_1}),f(e_{\epsilon_3})\right ).)\] Then
$T$ is an isomorphism of algebraic varieties, hence $\mathcal{V}$ is irreducible and
$\dim\mathcal{V} = 2(2n-(i+j))$.

As the coadjoint orbit $O_{\epsilon_i+\epsilon_j}(1)$ is an irreducible variety of dimension 
$2(2n-(i+j))$ also, we must have
 $\mathcal{V} = O_{\epsilon_i+\epsilon_j}(1)$.

\noindent \textbf{(ii)} In this case, we have 
\begin{equation*}
\begin{split}
S(\epsilon_i + \epsilon_j)  = & 
\bigcup_{k=i+1}^{j-1}\{\epsilon_i - \epsilon_k, \epsilon_k + \epsilon_j\} 
\cup \bigcup_{k=j+1}^{n}\{\epsilon_i - \epsilon_k, \epsilon_j + \epsilon_k\} \\ 
& \cup \bigcup_{k=j+1}^{n}\{\epsilon_i + \epsilon_k, \epsilon_j - \epsilon_k\} ,
\end{split}
\end{equation*}
and $\vert S(\epsilon_i + \epsilon_j) \vert =  2(2n-(i+j+1))$.
The set $S(\epsilon_i+\epsilon_j)$ can be written as the disjoint union of two subsets
$S_{(i)}$ and $S_{(j)}$ defined as follows:
\[ S_{(i)} = \bigcup_{k=i+1}^{j-1}\{\epsilon_i - \epsilon_k\} \cup \bigcup_{k=j+1}^n\{\epsilon_i-\epsilon_k\}\cup 
\bigcup_{k=j+1}^n\{\epsilon_i+\epsilon_k\}\] and
\[ S_{(j)} = \bigcup_{k=i+1}^{j-1}\{\epsilon_k+\epsilon_j\}\cup \bigcup_{k=j+1}^n\{\epsilon_j+\epsilon_k\}\cup 
\bigcup_{k=j+1}^n\{\epsilon_j-\epsilon_k\}.\]
(For example, if we have $\epsilon_1+\epsilon_3 \in \Phi_2^+(D_3)$, then
$S_{(1)}=\{\epsilon_1-\epsilon_2\}$ and $S_{(3)}=\{\epsilon_2+\epsilon_3\}$.)
As in Part \textbf{(i)}, for $\gamma \in S_{(i)}$, let $\gamma'$ be the unique element of $S_{(j)}$ such that
$\gamma + \gamma' = \epsilon_i + \epsilon_j$. Let $n(\gamma,\gamma')$ be the integer 
such that $[e_{\gamma},e_{\gamma'}] = n(\gamma,\gamma')e_{\epsilon_i+\epsilon_j}$
It follows from equations~(\ref{E:Drootvectors}) that $n(\gamma,\gamma') = \pm 1$.
 
 Let $\mathcal{V}$ be the variety  consisiting of all $f \in \gee^*$satisfying equations~(\ref{E:orbD2}).
 As in Part \textbf{(ii)}, we have that 
 $\mathcal{V} = O_{\epsilon_i+\epsilon_j}(1)$ . \qed

\bigskip

\noindent \textbf{Examples}

\noindent \textbf{(i)} Let $\gee = A_{3}^+$ and let 
$\alpha = \epsilon_1-\epsilon_4$. 
Then 
\[S(\alpha) = \{\epsilon_1-\epsilon_2,\epsilon_2-\epsilon_4,
\epsilon_1-\epsilon_3,\epsilon_3-\epsilon_4\},\] 
\[R(\alpha) = \{\epsilon_1-\epsilon_4, \epsilon_2-\epsilon_3\}.\]
So $O_{\alpha}(1)$ consists of
all $f$ which satisfy the equations
\[ f(e_{\epsilon_1-\epsilon_4}) = 1, \
f(e_{\epsilon_2-\epsilon_3})=f(e_{\epsilon_1-\epsilon_3})f(e_{\epsilon_2-\epsilon_4})\]
and $f$ takes arbitrary values on $e_{\beta}$ for all
$\beta \in S(\alpha)$.

\noindent \textbf{(ii)} Let $\gee = B_3^+$. Consider the roots $\epsilon_1-\epsilon_3,\epsilon_1,\epsilon_1+\epsilon_3 \in \Phi^+(B_3)$. Then we have
\[S(\epsilon_1-\epsilon_3)=\{\epsilon_1-\epsilon_2,\epsilon_2-\epsilon_3\},\]
\[R(\epsilon_1-\epsilon_3)=\{\epsilon_1-\epsilon_3,\epsilon_1,\epsilon_2,\epsilon_3,\epsilon_1+\epsilon_3,\epsilon_2+\epsilon_3,\epsilon_1
+\epsilon_2\},\]
\[S(\epsilon_1)=\{\epsilon_1-\epsilon_2,\epsilon_2,\epsilon_1-\epsilon_3,\epsilon_3\},\]
\[R(\epsilon_1)=\{\epsilon_2-\epsilon_3, \epsilon_2+\epsilon_3,\epsilon_1+\epsilon_3,
\epsilon_1+\epsilon_2\},\]
\[S(\epsilon_1+\epsilon_3)=\{\epsilon_1-\epsilon_2,\epsilon_2+\epsilon_3,\epsilon_1,\epsilon_3\},\]
\[R(\epsilon_1+\epsilon_3)=\{\epsilon_1+\epsilon_3,\epsilon_1-\epsilon_3,\epsilon_2-\epsilon_3,\epsilon_2, \epsilon_1+\epsilon_2\}.\]
Then $O_{\epsilon_1-\epsilon_3}(1)$ consists of all $f \in \gee^*$ that satisfy the equations
\[ f(e_{\epsilon_1-\epsilon_3})=1,\ f(e_{\gamma})=0\ \text{for all other}\ \gamma \in R(\epsilon_1-\epsilon_3)\] and 
$f$ takes arbitrary values on $e_{\beta}$ for all $\beta \in S(\epsilon_1-\epsilon_3)$.

\noindent The elementary orbit $O_{\epsilon_1}(1)$ consists of all $f \in \gee^*$ that satisfy the equations
\[ f(e_{\epsilon_1})=1, f(e_{\epsilon_2-\epsilon_3})=f(e_{\epsilon_1-\epsilon_3})f(e_{\epsilon_2}),\
f(e_{\gamma})=0\ \text{for all other}\  \gamma \in R(\epsilon_1),\] and $f$ takes arbitrary values on
$e_{\beta}$ for all $\beta \in S(\epsilon_1)$.

\noindent The elementary orbit $O_{\epsilon_1+\epsilon_3}(1)$ consists of all $f \in \gee^*$ that satisfy the equations
\[f(e_{\epsilon_1+\epsilon_3})=1,\ f(e_{\epsilon_1-\epsilon_3})=-\frac 12 f(e_{\epsilon_1})^2,\
f(e_{\epsilon_2-\epsilon_3})=f(e_{\epsilon_2+\epsilon_3})(-\frac 12 f(e_{\epsilon_2})^2),\]
\[f(e_{\epsilon_2})=f(e_{\epsilon_1})f(e_{\epsilon_2+\epsilon_3}),\ f(e_{\epsilon_1+\epsilon_2})=0,\]
and $f$ takes arbitrary values on $e_{\beta}$ for all $\beta \in S(\epsilon_1+\epsilon_3)$.

\noindent \textbf{(iii)} Let $\gee = D_3^+$. Consider the roots $\epsilon_1-\epsilon_3,\epsilon_1+\epsilon_3 \in
\Phi^+(D_3)$. Then we have
\[S(\epsilon_1-\epsilon_3)=\{\epsilon_1-\epsilon_2,\epsilon_2-\epsilon_3\},\]
\[R(\epsilon_1-\epsilon_3)=\{\epsilon_1-\epsilon_3,\epsilon_1+\epsilon_3,\epsilon_2+\epsilon_3,\epsilon_1+\epsilon_2\},\]
\[S(\epsilon_1+\epsilon_3)=\{\epsilon_1-\epsilon_2,\epsilon_2+\epsilon_3\},\] and
\[R(\epsilon_1+\epsilon_3)=\{\epsilon_1\pm \epsilon_3,\epsilon_2-\epsilon_3,\epsilon_1+\epsilon_2\}.\]
Then  $O_{\epsilon_1-\epsilon_3}(1)$ consists of all $f \in \gee^*$ that satisfy the equations
\[ f(e_{\epsilon_1-\epsilon_3})=1,\ f(e_{\gamma})=0\ \text{for all other}\ \gamma \in R(\epsilon_1-\epsilon_3),\] and 
$f$ takes arbitrary values on $e_{\beta}$ for all $\beta \in S(\epsilon_1-\epsilon_3)$.

\noindent The elementary orbit $O_{\epsilon_1+\epsilon_3}(1)$ consists of all $f \in \gee^*$ that satisfy the equations
\[ f(e_{\epsilon_1+\epsilon_3})=1,\ f(e_{\epsilon_1-\epsilon_3})= 1,\ f(e_{\epsilon_2-\epsilon_3})=
f(e_{\epsilon_2+\epsilon_3}),\ f(e_{\epsilon_1+\epsilon_2})=0,\]
$f(e_{\beta})$ is arbitrary for all $\beta \in S(\epsilon_1+\epsilon_3)$.
 
\section{Basic Sums of Elementary Orbits for $A_{n-1}^+$}\label{dimorbtypeA}
Here, we consider the set of positive roots $\apos$ and the corresponding nilpotent Lie algebra
$\gee = A_{n-1}^+$ of strictly upper triangular matrices. Then $G = \exp(\gee)$ is the group of all $n \times n$ unitriangular 
matrices with entries in $\C$. 
For $f \in \gee^*$, we define Supp$(f)$, the \emph{support} of $f$, as follows:
\[ \textrm{Supp}(f) = \{\alpha \in \apos \mid f(e_{\alpha}) \neq 0\}. \]
A subset $\dee \subset \apos$ is called a \emph{basic subset} if 
$\alpha - \beta \notin \apos$ for any $\alpha, \beta \in \dee$.
For example, $\dee = \{\epsilon_1-\epsilon_3,\epsilon_2-\epsilon_5,\epsilon_3-\epsilon_4\}$ is a basic
 subset of $\Phi^+(A_5)$. In particular, the empty set is a basic subset of $\apos$.

Let $\Cash$ denote the set of non-zero complex numbers.
Given a non-empty basic subset $\dee \subset \apos$ and a map $\phi : \dee \rightarrow \Cash$, we define the 
\emph{basic sum} $O_{\dee}(\phi)$ to be the set
\begin{equation}\label{E:basicsum}
 O_{\dee}(\phi) = \sum_{\alpha \in \dee}O_{\alpha}(\phi(\alpha)) \subset \gee^*.
  \end{equation}
If $\dee$ is empty, we may consider the empty function $\phi: \dee \rightarrow \Cash$, and in this
case we define $O_{\dee}(\phi) = 0$.
If $\dee$ is a basic subset of $\apos$, we define by $S(\dee)$
the subset
\begin{equation}\label{E:deesing} 
S(\dee) = \bigcup_{\alpha \in \dee}S(\alpha)
 \end{equation}
of $\apos$. A root $\alpha \in \apos$ will be called a $\dee$-\emph{singular root} if $\alpha \in S(\dee)$. In particular, if $\dee = 
\emptyset$ then we 
have $S(\dee) = \emptyset$.

We recall some important results:

\begin{pro}{\rm {\cite{Andre:ja95}}}\label{P:andre}
\begin{enumerate}
\item $O_{\dee}(\phi)$ is an irreducible subvariety of $\gee^*$, of
  dimension $s(\dee) := \vert S(\dee) \vert$.
\item For any $f \in \gee^*$ there exists a unique basic subset $\dee$
  of $\apos$ and a unique map 
$\phi: \dee \rightarrow \Cash$ such that $f \in O_{\dee}(\phi)$.
\end{enumerate}
\end{pro}

\bigskip

\noindent \textbf{Proof.}
\begin{enumerate}
\item Irreducibility of $O_{\dee}(\phi)$ is shown in {\cite[Sec.~2,
    Cor.~1]{Andre:ja95}} and its dimension is given by {\cite[Sec.~3,
    Thm.~2]{Andre:ja95}}.
\item The fact that $f$ is contained in $O_{\dee}(\phi)$ for a unique
  $\dee, \phi$ is proved in {\cite[Sec.~2, Prop.~3 \& Prop.~4]{Andre:ja95}}.
\end{enumerate}
   
\section{Homogeneous Basic Subvarieties}
The definitions and results in this subsection are taken from
{\cite[Sec.4]{Andre:ja95}}.

For any basic subset $\dee \subset \apos$ and any map $\phi \colon
\dee \rightarrow \Cash$, the basic sum $O_{\dee}(\phi)$ is
invariant under the coadjoint action of $G$, hence it is a union of coadjoint orbits.
What are the pairs $(\dee,\phi)$ for which 
$O_{\dee}(\phi)$ is a single coadjoint orbit? That depends
on the geometric configuration of the basic subset $\dee$ of $\apos$.
We start with the following definition:

A subset $C$ of $\apos$ is called a $\emph{chain}$ if
\[ C = \{\epsilon_{i_1}-\epsilon_{i_2},\epsilon_{i_2}-\epsilon_{i_3},\dots,\epsilon_{i_{r-1}}-\epsilon_{i_r}\}. \]
The cardinality $\vert C \vert$ is referred to as the $\emph{length}$
of the chain $C$. It is clear that a chain $C$ is a basic subset of $\apos$.
Now we define the $\dee$-\emph{derived roots} for any basic subset $\dee$ of 
$\apos$ as follows.

Let $C = \{\epsilon_{i_1}-\epsilon_{i_2},\epsilon_{i_2}-\epsilon_{i_3},\dots,\epsilon_{i_{r-1}}-\epsilon_{i_r}\}$, and 
$C' = \{\epsilon_{j_1}-\epsilon_{j_2},\epsilon_{j_2}-\epsilon_{j_3},\dots,\epsilon_{j_{s-1}}-\epsilon_{j_s}\}$ 
be two chains in $\dee$. Then the
pair $(C,C')$ will be called a \emph{special pair of chains} (with respect 
to $\dee$) if the following conditions are satisfied:
\begin{description}
\item[(i)]
$C,C'$ have the same length, i.e. $r = s$.
\item[(ii)]
$C,C'$ intertwine, i.e. $i_1 < j_1 < i_2 < j_2 < \dots < i_r < j_r$.
\item[(iii)]
If there exists $j_0, 1 \leq j_0 < j_1$, such that $\epsilon_{j_0}-\epsilon_{j_1} \in \dee$,
then $i_1 < j_0$.
\item[(iv)]
If there exists $i_{r+1}, i_r < i_{r+1} \leq n$, such that 
$\epsilon_{i_r}-\epsilon_{i_{r+1}} \in \dee$, then $i_{r+1} < j_r$. 
\end{description}

Then the root $\epsilon_{i_1}-\epsilon_{j_1}$ is called the $(C,C')$-derived root. In general,
a root $\epsilon_i-\epsilon_j \in \apos$ will be called a $\dee$-\emph{derived root} if there 
exists a special pair of chains $(C,C')$ in $\dee$ such that $\epsilon_i-\epsilon_j$ is
the $(C,C')$-derived root. The set of all $\dee$-derived roots is referred
to as the \emph{derived set of} $\dee$ and it is denoted by $\dee'$. It is clear
that $\dee' \subset S(\dee)$. For example, let $\dee$ be the basic set 
$\{\epsilon_1-\epsilon_3,\epsilon_3-\epsilon_5,\epsilon_2-\epsilon_4,\epsilon_4-\epsilon_6\} \subset
\Phi^+(A_5)$.
Then $\dee$ contains a special pair of chains $(C,C')$ where
\[ C = \{\epsilon_1-\epsilon_3,\epsilon_3-\epsilon_5\} \] and
\[ C'=\{\epsilon_2-\epsilon_4,\epsilon_4-\epsilon_6\}. \]
We see that $\epsilon_1-\epsilon_2$ is the $(C,C')$-derived root, so $\epsilon_1-\epsilon_2 \in \dee'$.

Below we recall the necessary and sufficient condition for $O_{\dee}(\phi)$ to be
a single coadjoint orbit.
\begin{thm}{ \rm ({\cite[Sec.~4, Thm.~3]{Andre:ja95}})}
Let $\dee$ be a basic subset of $\apos$ and let $\phi \colon \dee \rightarrow \Cash$ be a map.
Then $O_{\dee}(\phi)$ is a single coadjoint orbit
if and only if the derived set $\dee'$ of $\dee$ is empty.
\end{thm}
(In this case, $O_{\dee}(\phi)$ is the coadjoint orbit of 
$f = \sum_{\alpha \in \dee}\phi(\alpha)e_{\alpha}^*$.)

\section{Dimensions of Coadjoint Orbits for $\gee = A_{n-1}^+$}
In this section, we describe the properties of basic subsets that contain two-dimensional coadjoint
orbits and then determine all the possible dimensions of coadjoint orbits in $\gee^*$, for
$\gee = A_{n-1}^+$.
First, we prove the following lemma.
\begin{lem}\label{Alem1}
Let $\gee = A_{n-1}^+$. If $f \in \gee^*$ and the coadjoint orbit containing $f$ is
two-dimensional, then $\textrm{Supp}(f) \subset
\{\epsilon_i - \epsilon_j \in \apos \mid 1 \leq j-i \leq 2 \}$.
\end{lem}

\noindent \textbf{Proof.}
Suppose $\epsilon_r-\epsilon_s \in \textrm{Supp}(f)$ and $s-r > 2$.
Then 
$\vert S(\epsilon_r - \epsilon_s) \vert = 2(s-r-1) \geq 4$.
Thus, $\gee^f  \subseteq span\{e_{\beta} \mid \beta \in \apos, \beta \notin S(\epsilon_r-\epsilon_s) \}$.
So, $\dim\Omega_f = \textrm{codim}(\gee^f) \geq 4$, which is a
contradiction of the hypothesis. \qed

\bigskip

Using this lemma, we can show the following:
\begin{thm}\label{Ares1}
Assume $\gee=A_{n-1}^+$. Let $f \in \gee^*$, and let $(\dee,\phi)$ be the unique pair such
that $f \in O_{\dee}(\phi)$ where $O_{\dee}(\phi)$ is as in (\ref{E:basicsum}). Let $\Omega_f$ denote the coadjoint orbit that
contains $f$. If $\dim(\Omega_f) = 2$,
 then $s(\dee) =\vert S(\dee) \vert = 2$ or $3$ where $S(\dee)$ is as in (\ref{E:deesing}).
\end{thm}
\noindent \textbf{Proof.}
Since $\dim\Omega_f = 2$, therefore $f \neq 0$ and  $s(\dee) \geq 1$.

Suppose $s(\dee) \geq 4$. Then there are $\epsilon_{i_1}-\epsilon_{j_1}, \epsilon_{i_2}-\epsilon_{j_2} \in \dee$ such
that
$j_1 - i_1 = 2 = j_2 - i_2$ and $i_2 - i_1 \geq 1$.
Again, $\{e_{\beta} \mid \beta \in S(\epsilon_{i_1}-\epsilon_{j_1}) \cup S(\epsilon_{i_2}-\epsilon_{j_2})\}$ 
is not contained in $\gee^f$, so $\dim\Omega_f \geq 4$, which is a contradiction.
\qed

\begin{lem}\label{Alem2}
There exist basic
subsets $\dee_m$ of $\apos$ such that the derived set of $\dee_m$ is empty, 
and $s(\dee_m) = 2m$ where 
\begin{description}
\item[(i)] 
$m=0,1,2,3,\dots,\frac{1}{4}(n-2)(n)$ if $n$ is even,
\item[(ii)]
$m = 0,1,2,3,\dots,\frac{1}{4}(n-1)^2$ if $n$ is odd.
\end{description}
\end{lem}

\noindent \textbf{Proof.}
If $n = 2$ then $\Phi^+(A_1)=\{\epsilon_1-\epsilon_2\}$. The only non-empty basic set $\dee_0 = \{\epsilon_1-\epsilon_2\}$ 
has no derived roots and $s(\dee_0) = 0 = \frac{1}{4}(2-2)(2)$, so this case
is trivial.

We will use induction on $n \geq 3$.
We start with the cases $n = 3$ and $n = 4$.

\bigskip
$\Phi^+(A_2) = \{\epsilon_1-\epsilon_2, \epsilon_2-\epsilon_3, \epsilon_1-\epsilon_3\}$:

\smallskip
\noindent The basic subset $\dee_0 = \{\epsilon_2-\epsilon_3\}$ has no derived roots and $s(\dee_0) = 0$.
The basic subset $\dee_1 = \{\epsilon_1-\epsilon_3\}$ has no derived roots and 
$s(\dee_1) = 2  = \frac{1}{2}(3-1)^2$.
We observe that there are no basic sets $\dee$ with $s(\dee)$ higher than 2.

\bigskip
$\Phi^+(A_3) = \{\epsilon_1-\epsilon_2, \epsilon_2-\epsilon_3,\epsilon_3-\epsilon_4, \epsilon_1-\epsilon_3, \epsilon_2-\epsilon_4,
 \epsilon_1-\epsilon_4\}$:
 
\smallskip
\noindent The basic subset $\dee_0 = \{\epsilon_2-\epsilon_3\}$ has no derived roots and $s(\dee_0) = 0$.
The basic subset $\dee_1 = \{\epsilon_2-\epsilon_4\}$ has no derived roots and $s(\dee_1) = 2$.
The basic subset $\dee_2 = \{\epsilon_1-\epsilon_4\}$ has no derived roots and $s(\dee_2) = 4
= \frac{1}{2}(4-2)(4)$. We observe that there are no basic sets $\dee$ with
$s(\dee)$ higher than 4.

\bigskip
Now let $n > 4$ be even.
Assume that \textbf{(i)} is true for all even numbers bigger than 3
and less than 
$n$. Let 
\[ \dee_0 = \{\epsilon_{\frac{n}{2}}-\epsilon_{\frac{n}{2}+1}\}, \quad
\dee_1 = \{\epsilon_{\frac{n}{2}}-\epsilon_{\frac{n}{2} + 2}\}, \]
\[\dee_2 = \{\epsilon_{\frac{n}{2} - 1}-\epsilon_{\frac{n}{2} + 2}\},
\dots,
\dee_{n-2} = \{\epsilon_1 - \epsilon_n\}.\]
(For example, if $n=6$, then $\dee_0 =\{\epsilon_3-\epsilon_4\},\dee_1=\{\epsilon_3-\epsilon_5\},
\dee_2=\{\epsilon_2-\epsilon_5\}, \dee_3=\{\epsilon_2-\epsilon_6\},\dee_4=\{\epsilon_1-\epsilon_6\}$.)
 
These are basic subsets of $\apos$
whose derived sets are empty and we see that 
$s(\dee_r) = 2r$ for $r = 0,1,2,\dots,n-2$.
The set $ \{\epsilon_i-\epsilon_j \mid  2 \leq i \leq j \leq n-1\} \subset \apos$ can be identified with the
set of roots $\Phi^+(A_{n-3})$. But $n-2$ is an even number smaller than $n$,
so \textbf{(i)} holds for $\Phi^+(A_{n-3})$. Thus, there exist basic subsets 
$\hat{\dee}_m$ of
$\Phi^+(A_{n-3})$ whose derived sets are empty, such that 
$s(\hat{\dee}_m)=2m$, where
$m=0,1,2,3,\dots,\frac{1}{4}(n-4)(n-2)$. If $\hat{\dee}_m$ is a basic subset of
$\Phi^+(A_{n-3})$ which has no derived roots, then $\dee_m =\hat{\dee}_m \cup \{\epsilon_1-\epsilon_n\}$ is
a basic subset of $\apos$ which has no derived roots, and also
$s(\dee_m) = 2(n-2) + s(\hat{\dee}_m)$. Thus, there exist basic subsets $\dee_m$ of 
$\apos$
whose derived sets are empty such that $s(\dee_m) = 2m$, where
$m=0,1,2,3,\dots,n-2,(n-2) + 1, (n-2) + 2, \dots, (n-2) + 
\frac{1}{4}(n-4)(n-2)$. But $(n-2) + \frac{1}{4}(n-4)(n-2) = 
\frac{1}{4}(n-2)(n)$, so the result holds in the even case.

Now let $n > 3$ be odd, and
assume that \textbf{(ii)} is true for all odd numbers less than $n$.
Let 
\begin{gather}
 \dee_0 = \{\epsilon_{\frac{n+1}{2}}-\epsilon_{\frac{n+1}{2}+1}\},\quad
\dee_1 = \{\epsilon_{\frac{n+1}{2} - 1}-\epsilon_{\frac{n+1}{2}+1} \}, \notag \\
\dee_2 = \{\epsilon_{\frac{n+1}{2} - 1}-\epsilon_{\frac{n+1}{2} + 2} \}, \dots, 
\dee_{n-2} = \{\epsilon_1-\epsilon_n\}. \notag 
\end{gather}
(For example, if $n=5$ then $\dee_0=\{\epsilon_3-\epsilon_4\},\dee_1=\{\epsilon_2-\epsilon_4\},
\dee_2=\{\epsilon_2-\epsilon_5\}, \dee_3=\{\epsilon_1-\epsilon_5\}$.)

 These are basic subsets of $\apos$
which have no derived roots and we see that 
$s(\dee_r) = 2r$ for $r = 0,1,2,\dots,n-2$.
The set
$\{\epsilon_i-\epsilon_j \mid 2 \leq i \leq  j \leq n-1 \} \subset \apos$
can be identified with $\Phi^+(A_{n-3})$. Since $n-2$ is an odd number 
less than $n$, \textbf{(ii)} is true for $\Phi^+(A_{n-3})$. Therefore by the induction hypothesis, there exist basic
subsets $\hat{\dee}_m$ of $\Phi^+(A_{n-3})$ which have no derived roots, such that
$s(\hat{\dee}_m)=2m$, where $m=0,1,2,3,\dots,\frac{1}{4}(n-3)^2$. Now, if
$\hat{\dee}_m$ is a basic subset of $\Phi^+(A_{n-1})$ which has no derived roots,
then $ \dee_m = \hat{\dee}_m  \cup \{\epsilon_1-\epsilon_n\}$ is a basic subset of $\apos$ which 
has
no derived roots. Also, $s(\dee_m) = 2(n-2) + s(\hat{\dee}_m)$.
Thus, there exist basic subsets $\dee_m$ of $\apos$ which have no derived
roots such that $s(\dee_m)=2m$,where $m = 0,1,2,3,\dots,n-2,(n-2) + 1,
(n-2) + 2, \dots, 
(n-2) + \frac{1}{4}(n-3)^2$.
As $(n-2) + \frac{1}{4}(n-3)^2 = \frac{1}{4}(n-1)^2$, this gives the desired conclusion. \qed

\begin{thm}\label{orbres2}
There exist elements $ f_m \in (A_{n-1}^+)^*$ such that 
$\dim \Omega_{f_m}=2m$, where
\begin{description}
\item[(i)] 
$m=0,1,2,3,\dots,\frac{1}{4}(n-2)(n)$ if $n$ is even,
\item[(ii)]
$m=0,1,2,3,\dots,\frac{1}{4}(n-1)^2$ if $n$ is odd.
\end{description}
\end{thm}

\noindent \textbf{Proof.}
From Lemma~\ref{Alem2} and {\cite[Sec.~4, Thm.~3]{Andre:ja95}}, we see that
there exist basic subsets $\dee_m$ of $\apos$ such that for any
map $\phi_m \colon \dee_m \rightarrow \Cash$, the variety
$O_{\dee_m}(\phi_m)$ is a single coadjoint orbit containing 
$f_m = \sum_{\alpha \in \dee_m}\phi_m(\alpha)e_{\alpha}^*$. Since 
$\dim O_{\dee_m}(\phi_m) = s(\dee_m)$, we have $\dim \Omega_{f_m}=2m$, where
\begin{description}
\item[(i)]
$m = 0,1,2,3,\dots,\frac{1}{4}(n-2)n$ if $n$ is even,
\item[(ii)]
$m = 0,1,2,3,\dots,\frac{1}{4}(n-1)^2$ if $n$ is odd.
\end{description}
\qed

\begin{thm}\label{orbres3}
If $\gee = A_{n-1}^+$, then for any $f \in \gee^*$,
\[ \dim \Omega_f \leq
\begin{cases}
\frac{1}{2}(n-2)n &\text{if $n$ is even} \notag \\
\frac{1}{2}(n-1)^2 &\text{if $n$ is odd} \notag \\
\end{cases} \]
\end{thm}

\noindent \textbf{Proof.} Let $(\dee,\phi)$ be the unique pair such
that $f \in O_{\dee}(\phi)$. Then $\Omega_f \subset O_{\dee}(\phi)$,
therefore $\dim \Omega_f \leq s(\dee)$, by Proposition~\ref{P:andre}.

If $n$ is odd, then the basic subset $\dee =
\{\epsilon_1-\epsilon_n,\epsilon_2-
\epsilon_{n-1},\dots,\epsilon_{\frac{n+1}{2}-1}-\epsilon_{\frac{n+1}{2}+1}\}$
is the basic subset with the largest number of $\dee$-singular roots
(i.e. with the largest value of $s(\dee)$), and $s(\dee) = \frac 12(n-2)n$.

If $n$ is even, then the basic subset $\dee =
\{\epsilon_1-\epsilon_n,\epsilon_2-\epsilon_{n-1},
\dots,\epsilon_{\frac n2}-\epsilon_{\frac n2 +1}\}$
is the basic subset with the largest number of $\dee$-singular roots
(i.e. with the largest value of $s(\dee)$), and $s(\dee) = \frac{1}{2}(n-1)^2$.

Hence the theorem follows.\qed

\bigskip
\begin{cor}
Let $\gee = A_{n-1}^+$ and let $U$ be the universal enveloping algebra of $\gee$.
Then, for any primitive ideal $I$ of $U$ we have $U/I \simeq \mathcal{A}_m$ where
\begin{description}
\item[(i)] 
$m=0,1,2,3,\dots,\frac{1}{4}(n-2)n$ if $n$ is even,
\item[(ii)]
$m=0,1,2,3,\dots,\frac{1}{4}(n-1)^2$ if $n$ is odd;
\end{description}
and $\mathcal{A}_m$ is the $m$-th Weyl algebra.
\end{cor}

\noindent \textbf{Proof.} This is clear from Theorem~\ref{orbres2},
Theorem~\ref{orbres3}, and the remarks at the beginning of 
Section~\ref{S:elemcoadorb}. \qed

\end{document}